\documentclass{amsart}
\usepackage[OT4]{fontenc}
\usepackage{latexsym,color}
\usepackage{amsmath,amsthm,amsfonts,amssymb,graphicx}

\title{The asymptotic dimension of a curve graph is finite}
  \author{Gregory Bell}
  \address{Mathematics and Statistics, UNC Greensboro, Greensboro, NC 27402, USA}
  \email{gcbell@uncg.edu}
  \author{Koji Fujiwara}
  \address{Mathematical Institute, Tohoku University, Sendai, 980-8578 Japan}
  \email{fujiwara@math.tohoku.ac.jp}
  \date{November 1, 2006}

\newtheorem{thm}{Theorem}         
\newtheorem{prop}{Proposition} 
\newtheorem{cor}{Corollary}
\newtheorem{lemma}{Lemma}

{}

\theoremstyle{remark}

\newtheorem{question}{Question}
\newtheorem*{examples*}{Examples}
\newtheorem*{definition*}{Definition}
\newtheorem*{remark*}{Remark}
\newtheorem*{example*}{Example}
\newtheorem*{cor1}{Corollary 1}
\newtheorem*{cor3}{Corollary 3}
\newtheorem*{remarks*}{Remarks}
\DeclareMathOperator{\as}{asdim} \DeclareMathOperator{\diam}{diam}
\DeclareMathOperator{\fin}{fin} \DeclareMathOperator{\supp}{supp}
\DeclareMathOperator{\Mod}{Mod} \newcommand{\N}{\mathbb{N}}
\newcommand{\G}{\Gamma} \newcommand{\sU}{\mathcal{U}}
\newcommand{\sV}{\mathcal{V}} \newcommand{\sL}{\mathcal{L}}
\newcommand{\I}{\mathcal{I}}\newcommand{\R}{\mathbb{R}}
\newcommand{\sW}{\mathcal{W}}\newcommand{\sA}{\mathcal{A}}
\newcommand{\Z}{\mathbb{Z}} \DeclareMathOperator{\vcd}{vcd}

\dedicatory{Dedicated to Professor Yukio Matsumoto on his 60th
birthday.}

\subjclass[2000]{57M99 (primary), 20F69 (secondary)}

\begin{document}
\maketitle
\begin{abstract}
We find an upper bound for the asymptotic dimension of a hyperbolic
metric space with a set of geodesics satisfying a certain
boundedness condition studied by Bowditch. The primary example is a
collection of tight geodesics on the curve graph of a compact
orientable surface. We use this to conclude that a curve graph has
finite asymptotic dimension. It follows then that a curve graph has
property $A_1$. We also compute the asymptotic dimension of mapping
class groups of orientable surfaces with genus $\le
2$.\end{abstract}

\section{Introduction}

The asymptotic dimension of a metric space was defined by Gromov
\cite[page 29]{Gr93} as the large-scale analog of Lebesgue covering
dimension. Gromov originally denoted asymptotic dimension $\as_+,$
but the notation $\as$ has become more standard.  Define $\as X\le
n$ if for every (large) $r>0$ there is a cover of $X$ by uniformly
bounded subsets of $X$ so that no ball of radius $r$ in $X$ meets
more than $n+1$ elements of the cover, (see section 2). Asymptotic
dimension is a coarse invariant (see section 2), so in particular it
is a quasi-isometry invariant.

The notion of a $\delta$-hyperbolic space is due to Gromov
\cite{Gr87} and has been studied extensively. Gromov remarked
\cite[page 31]{Gr93} that finitely generated hyperbolic groups have
finite asymptotic dimension. Recently Roe \cite{Ro05} improved on
this result, showing that a hyperbolic geodesic metric space with
bounded growth (see section 2) has finite asymptotic dimension. This
also follows from a result of Bonk and Schramm \cite{BS}.

Below we consider the asymptotic dimension of hyperbolic geodesic
graphs which need not have bounded growth, but satisfy a boundedness
condition in terms of geodesics that we call \emph{property B}.
Although the property is stated locally, it is indeed a global
condition. A simple example of a space with property B that does not
have bounded growth is an infinite-valence tree. In Theorem
\ref{main} we prove that hyperbolic spaces with property B have
finite asymptotic dimension.

Let $S=S_{g,p}$ denote a compact, orientable surface of genus $g$
with $p$ punctures (for technical reasons, we require $3g-3+p>1$).
The \emph{complex of curves}, $C(S)$ was defined by Harvey
\cite{harvey} in his study of Teichm\"uller spaces of Riemannian
manifolds. It is a simplicial complex whose vertices are isotopy
classes of essential curves on $S$; a collection of vertices span a
simplex when they can be realized simultaneously by disjoint curves
on $S$. The \emph{curve graph} is the $1$-skeleton $C^{(1)}(S)$ and
is denoted $X(S)$. Masur and Minsky \cite{minsky-masur:cc1} showed
that $X(S)$ is a hyperbolic graph, (see Theorem \ref{thm2}).
Although $X(S)$ does not have bounded growth, Bowditch showed that
it does satisfy our boundedness condition, \cite{bowditch:tight}.
Therefore we can apply Theorem \ref{main} to $X(S)$ to obtain the
main result of the first part of the paper:

\begin{cor1} Let $S$ be an orientable surface with genus $g$ and $p$
punctures such that $3g-3+p>1$. Then $\as X(S)<\infty$.
\end{cor1}

A metric space is said to be \emph{proper} if closed, bounded sets
are compact. Although it is interesting to study asymptotic
dimension for its own sake, much of the original interest in
asymptotic dimension arose in connection with the Novikov
conjecture. G. Yu \cite{Yu1} showed that the coarse Baum-Connes
conjecture holds for proper metric spaces with finite asymptotic
dimension. Applying the so-called descent principle, he concludes
that the Novikov conjecture holds for finitely generated groups with
finite classifying space and finite asymptotic dimension. This
result has since been improved upon by Yu \cite{Yu2} and Kasparov-Yu
\cite{KY}.

In \cite{Yu2}, Yu defines a property of discrete metric spaces
called ``property $A$" guaranteeing the existence of a large-scale
(or uniform) embedding into Hilbert space. The existence of this
embedding implies the coarse Baum-Connes conjecture for discrete
spaces with bounded geometry. As a consequence, finitely generated
groups with this property satisfy the Novikov conjecture. Yu's
property $A$ is a coarse invariant \cite{Yu2,Tu}. Many of the
examples we consider in this paper are not proper and moreover, they
are not even coarsely equivalent to proper metric spaces, so the
coarse Baum-Connes conjecture does not apply, (see Proposition
\ref{curvegraphnotproper}). 

In the fifth section of this paper we show that spaces with our
boundedness condition have property $A$ in the sense of Tu
\cite{Tu}, which we call property $A_1$. This property reduces to
Yu's property $A$ when the spaces are discrete with bounded
geometry. In particular, we prove:

\begin{cor3} Let $S$ be an orientable surface with genus $g$ and $p$
punctures so that $3g-3+p>1$. Then $X(S)$ has property $A_1$.
\end{cor3}

This follows from Corollary 1 and the proof of a theorem of
Higson-Roe \cite[Lemma 4.3]{HR}.

In the sixth section of the paper, we consider the mapping class
group of $S$, denoted $\Mod(S)$. By combining a recent result of
Dranishnikov \cite{Dr05} with some basic facts about the mapping
class groups of surfaces with $g\le 2$ we can show that
$\as\Mod(S)=\vcd\Mod(S)$ when $g\le 2$ and $3g-3+p>1$. As a
corollary, we also find upper bounds for the asymptotic dimension of
braid groups, (Corollary \ref{braid}) and find an exact formula for
the asymptotic dimension of some Artin groups, (Corollary
\ref{artin-cor}).

We conclude the paper with a list of open questions about curve
complexes and mapping class groups.

\textbf{Acknowledgements.} We wish to thank A.~Dranishnikov and
G.~Yu for many helpful remarks and J.~Behrstock, B.~Bowditch,
D.~Margalit, S.~Schleimer and M.~Kapovich for very helpful comments
on a preliminary version of the paper. The second author would like
to thank K.~Bromberg for a helpful discussion and the Max Planck
Institute for Mathematics in Bonn for their hospitality. Finally,
the authors would like to thank organizers of the International
Conference on Geometric Topology in B\k{e}dlewo where collaboration
for this paper began.

We also wish to express our very sincere gratitude to the anonymous
referee for many helpful comments and a careful reading of the
manuscript.

\section{Preliminaries}

There are several equivalent definitions of asymptotic dimension of
a metric space (see the survey article \cite{BD4} or the book
\cite{Ro03}). In this paper we will use only the following
definition, (mentioned in the introduction) in terms of
$r$-multiplicity of a cover. In a metric space $X$, let $N(x;r)$
denote the $r$-ball centered at $x$ in $X.$ Define $\as X\le n$ if
for every (large) $r\in\mathbb{N}$ there exists a cover $\{U_i\}_i$
of $X$ by uniformly bounded subsets of $X$ so that $\sharp\{i\mid
N(x;r)\cap U_i\neq\emptyset\}\le n+1$ for all $x\in X.$ The number
$\sup_{x\in X}\sharp\{i\mid N(x;r)\cap U_i\neq\emptyset\}\le n+1$ is
called the $r$-\emph{multiplicity} of the cover $\{U_i\}_i$.

A geodesic triangle in a metric space $X$ is said to be
$\delta$-{\em thin} if each of its sides is contained in the
$\delta$-neighborhood of the union of the other two sides.  A
geodesic space $X$ will be called $\delta$-{\em hyperbolic} if every
geodesic triangle in $X$ is $\delta$-thin.  By a {\em hyperbolic
metric space} we mean a geodesic metric space which is
$\delta$-hyperbolic for some $\delta\ge 0.$  We will also say that a
finitely generated group is {\em hyperbolic} if the Cayley graph
with respect to some choice of a finite generating set is
$\delta$-hyperbolic for some $\delta\ge 0.$ Observe that this
property is a quasi-isometry invariant -- although the value of
$\delta$ is not -- so it makes sense to speak of a finitely
generated group as hyperbolic without reference to a metric.

A subset $A$ of a metric space $X$ is $r$-\emph{discrete} if
$d(a,a')\ge r$ for all $a\neq a'$ in $A.$ The $r$-\emph{capacity} of
a set $Y\subset X$ is the maximal cardinality of an $r$-discrete set
in $Y$. A metric space $X$ has \emph{bounded geometry} if there is
an $r>0$ and a function $c:[0,\infty)\to [0,\infty)$ so that for any
$x\in X$ the $r$-capacity of $N(x;R)$ does not exceed $c(R).$ For a
discrete metric space, this simply means that the cardinality of any
$R$-ball is bounded by $c(R).$

A metric space $X$ has {\em bounded growth at some scale} if there
are constants $r$ and $R$ with $R> r
> 0,$ and $N\in\N$ such that every open ball of radius $R$ in $X$
can be covered by $N$ open balls of radius $r.$ A metric space with
bounded geometry has bounded growth: take $r$ from the definition of
bounded geometry, take $R>r$ and put $N=c(R).$ Clearly the Cayley
graph of a finitely generated group has bounded geometry, and so has
bounded growth. Thus, any finitely generated group has bounded
growth at some scale.

We mentioned that asymptotic dimension is a quasi-isometry
invariant, but we will need a stronger statement. It is also a
coarse invariant \cite{Ro03}.  Let $X$ and $Y$ be metric spaces. A
map $f:X\to Y$ is \emph{bornologous} if for every $R>0$ there is an
$S>0$ such that $d(f(x),f(x'))<S$ whenever $d(x,x')<R.$ The map $f$
is \emph{metrically proper} if $f^{-1}(B)$ has bounded diameter in
$X$ for each bounded subset $B\subset Y.$ A \emph{coarse map} is a
map $f:X\to Y$ that is bornologous and metrically proper. A coarse
map $f:X\to Y$ is a \emph{coarse equivalence} between $X$ and $Y$ if
there is a coarse map $g:Y\to X$ and a constant $K\ge 0$ so that
$d(fg(y),y)\le K$ and $d(gf(x),x)\le K$ for all $x\in X$ and $y\in
Y.$ It is worth mentioning that in the case that $X$ and $Y$ are
finitely generated groups -- or indeed any length metric spaces --
the notions of coarse equivalence and quasi-isometry coincide, see
\cite{Ro03}.

\section{Boundedness condition}

Let $\G$ be a $\delta$-hyperbolic graph. Suppose that $\sL$ is a set
of geodesics in $\G$ such that for any $a,b \in V(\G)$, there exists
$\gamma \in \sL$ connecting $a$ and $b$.

For $a,b \in V(\G)$, let $\sL(a,b)$ be the set of all geodesics in
$\sL$ connecting $a$ and $b$. We write $G(a,b)= \cup \sL(a,b)
\subset \G$. Given $A,B \subset V(\Gamma)$, let $\sL(A,B)=\cup
\{\sL(a,b) \mid a \in A, b \in B\}$ and $G(A,B)=\cup \sL(A,B)$. We
write $G(a,b;r)$ to mean the set $G(N(a;r),N(b;r))$.

In \cite{bowditch:tight}, Bowditch considered the following
boundedness property, which we call {\it property B}\,: there are
constants $\ell,k,D$ such that if $a,b \in V(\Gamma), r \in {\Bbb
N}$ and $c \in G(a,b)$ with $d(c, \{a,b \}) \ge r+\ell$, then
$G(a,b;r) \cap N(c;k)$ has at most $D$ elements. We also assume that
any geodesic connecting $N(a;r)$ to $N(b;r)$ must intersect
$N(c;k)$.

We remark that if this property is satisfied, then we can assume
that $k=2 \delta$ because $\G$ is $\delta$-hyperbolic.

We say that a space $X$ has {\em property B} if there is a
collection $\sL$ of geodesics on the space with the above
properties. 
%

\begin{thm}\label{main}
Let $\G$ be a $\delta$-hyperbolic graph. Suppose $\sL$ is a set of
geodesics in $\G$ such that any two vertices of $\G$ are joined by a
geodesic in $\sL$. If $\sL$ satisfies property B with constants
$\ell,k=2 \delta$ and $D$, then the asymptotic dimension of $\Gamma$
is at most $2D-1$.
\end{thm}
\begin{proof} Suppose $ r \in \N$ is given. We may assume that $\ell \in
\N$, and $\ell \ge 10 \delta$. Fix $x_0 \in \G$. For each $n \in
\N$, define
$$A_n=\{x \in \G\mid10(n-1)(r+\ell) \le d(x,x_0) \le 10n(r+\ell) \}.$$
Clearly $\cup _n A_n =\G$. Let
$$S_n=\{x \in \G\mid d(x,x_0)=10n(r+\ell) \}.$$
For each $n \ge 3$, we define subsets $\{U^n_i \}_i$ of $A_n$ such
that $\cup _i U^n_i =A_n$ as follows. Write out the elements of
$S_{n-2}$ as $\{s_i\}$. Define
$$U^n_i=\{ x \in A_n\mid ^{\exists}[x,x_0] \in \sL \, {\rm s.t.} \, s_i \in [x,x_0] \},$$
where $[x,x_0]$ is a geodesic from $x$ to $x_0$. If $n=1$ or $2$,
let $U^n_1=A_n$. Clearly, $\cup_i U^n_i =A_n$ for each $n$, so the
collection $\{U^n_i\}_{i,n}$ covers $\G$.

\noindent {\bf Claim:} (1) for any $i$ and for each $n$, $\diam
U^n_i \le 40(r+\ell)$; and
\\
(2) for an $r/2$-ball $V$ in $\G$ and for each $n$, there are at
most $D$ subsets $U^n_i$ such that $U^n_i \cap V \not= \emptyset$.

To see (1), suppose $x,y \in U^n_i$. If $n=1$ or $2$ this is clear,
so suppose $n \ge 3$. Then, $d(x,s_i) \le 20(r+\ell),$ and $d(y,s_i)
\le 20(r+\ell)$, so that $d(x,y) \le 40(r+\ell)$.

For (2), we observe that this is clear if $n=1$ or $2$, so suppose
$n \ge 3$. Suppose $U^n_i \cap V \not= \emptyset$ and $U^n_j \cap V
\not=\emptyset$. Choose $y_i \in U^n_i \cap V, y_j \in U^n_j \cap
V$. Then, $d(y_i,y_j) \le r$, so that $y_j \in N(y_i,r)$. Therefore,
$s_j \in G(y_i,x_0;r)$, because $s_j \in [x_0,y_j] \in \sL$. By the
$\delta$-hyperbolicity of $\G$, $s_j \in N(s_i;2 \delta)$, since
$d(y_i,s_i) \ge r+\ell$ and $\ell \ge 10 \delta$. We find $s_j \in
G(y_i,x_0;r) \cap N(s_i;2 \delta)$. Fixing $i$ and allowing $j$ to
vary we see, by property B (letting $x_0=a, y_i=b, s_i=c$), that
there are at most $D$ such $s_j$ satisfying this property since
$d(x_0,s_i) \ge r+\ell$. So we have the claim.

It follows that the collection $\{U^n_i \}_{i,n}$ is a uniformly
bounded cover with $r/2$-multiplicity $\le 2D.$  The last item is
clear from claim (2) because for any $x\in\G$ there can be at most
two sets of the form $A_n$ with $A_n \cap N(x;r/2) \not= \emptyset$.
Thus, $\as\G\le 2D-1$ as required. \end{proof}

\begin{examples*}
\begin{enumerate}
\item
Let $T$ be a tree, which is $0$-hyperbolic. Let $\sL$ be the set of
all geodesics in $T$. $\sL$ satisfies the property B with $D=1, k=0,
\ell=0$, so that the asymptotic dimension of $T$ is at most $1$.  We
should note that this was known to Gromov \cite[\textbf{1.E.},
Example (b)]{Gr93} and an explicit proof appears in Roe's book
\cite{Ro03}.
\item
Let $\G$ be the Cayley graph of a word-hyperbolic group $G$ with
respect to a finite generating set, $S$, with $S^{-1}=S$. Let
$|S|=s$. Let $\sL$ be the set of all geodesics in $\G$. Any two
points are joined by a geodesic. Suppose that $\G$ is
$\delta$-hyperbolic. For a point $a \in V(\G)$, let $D=|N(a;2
\delta)| \le s^{2 \delta}$. Note that $D$ does not depend on $a$
since $G$ acts transitively on $\G$ by isometries. Then, $\sL$
satisfies property B with $k=2 \delta, \ell=10 \delta$ and $D$.
Thus, we conclude that $\as\G\le 2s^{2\delta}-1.$
\end{enumerate}
\end{examples*}

\section{Curve graph}
Let $S=S_{g,p}$ be a compact orientable surface of genus $g$ with
$p$ punctures. The {\it curve complex} of $S$ was defined by Harvey
\cite{harvey} and has been successfully used in the study of mapping
class groups, for example in \cite{BF,Hare,I}. For a good overview
we recommend Ivanov \cite{Iv}. For our purposes, we will restrict
attention to the 1-skeleton of the curve complex, called the
\emph{curve graph}. The curve graph $X$ of $S$ is a graph whose
vertices are isotopy classes of essential, nonparallel,
nonperipheral, simple closed curves in $S$ where two distinct
vertices are joined by an edge if the corresponding curves can be
realized simultaneously by pairwise disjoint curves. In certain
sporadic cases $X$ as defined above is $0$-dimensional or empty.
This happens when there are no curve systems consisting of two
curves, i.e. $3g-3+p \le 1$, (see the discussion following the
Remark (ii) below). 
In the discussion that follows these cases are excluded. The mapping
class group $Mod(S)$ of $S$ (see Section 6) acts on $X$ by $f\cdot
a=f(a)$. H. Masur and Y. Minsky proved the following remarkable
result.
\begin{thm}[\cite{minsky-masur:cc1}]\label{thm2}
Let $S$ be a compact orientable surface of genus $g$ with $p$
punctures. Suppose $3g-3+p>1$. Then the curve graph $X(S)$ is
hyperbolic.
\end{thm}

It is known (see \cite{minsky-masur:cc1}) that $X$ is connected,
locally infinite, and its diameter is infinite. In general there are
infinitely many geodesics connecting $a,b\in X$; however,
Masur-Minsky found a set, $\sL$, of geodesics, called {\it tight
geodesics} such that $\sL(a,b)$ is not empty and finite for any $a,b
\in X$.

Bowditch showed the following.
\begin{thm} {\rm\cite[Theorem 1.2]{bowditch:tight}}\label{thm3}
Let $X$ denote the curve graph of a surface $S_{g,p}$ with
$3g-3+p>1.$ Let $\sL$ be the set of tight geodesics on $X.$ Then
$\sL$ satisfies property B for some $k,\ell$ and $D.$ 
\end{thm}

It follows from Theorem \ref{main} that the asymptotic dimension of
the curve graph $X$ is at most $2D-1$ so in particular it is finite.
\begin{cor}
Let $S_{g,p}$ be an orientable surface of genus $g$ with $p$
punctures
such that $3g-3+p > 1$. Then $\as X(S_{g,p})<\infty$.
\end{cor}

\begin{remarks*}
\begin{enumerate}
\item As we said, the curve graph of $S$ is the $1$-skeleton of the
curve complex of $S$. They are quasi-isometric, so they have same
asymptotic dimension.
\item After this paper was written, Fujiwara and Whyte \cite{FW} showed
that a hyperbolic geodesic space with asymptotic dimension $1$ is
quasi-isometric to a tree. Using results of Behrstock and Leininger,
Schleimer \cite{Schl} showed the following property for
$X=X(S_{g,1})$ with $g \ge 2$: for any $x \in X$ and $r \in {\Bbb
N}$, $X \backslash N(x;r)$ is connected. It follows that $X$ is not
coarsely equivalent to a tree (cf. Proposition
\ref{curvegraphnotproper}). This implies that $\as X(S_{g,1})>1$
($g\ge 2$), \cite{FW}. Thus there exist surfaces whose curve graphs
have asymptotic dimension greater than 1. It is unclear what happens
in other cases.
\end{enumerate}
\end{remarks*}

As mentioned above, when $3g-3+p\le 1$ the curve graph as defined
above is $0$-dimensional or empty. In these cases, one could rectify
the situation by declaring two vertices to be joined by an edge if
the corresponding curves can be realized with only one intersection
point. One such case of particular interest is the once-punctured
torus ($g=1,\,p=1$). By placing an edge between curves if their
intersection number is $1$, we obtain a Farey graph. It is the
$1$-skeleton of a planar $2$-complex whose dual tree is a binary
tree.

A Farey graph $X$ (or the Farey graph since they are all isomorphic)
can also be described as follows (see, for example \cite[\S
3]{Min96}, which provides a nice picture of $X$). In the rational
projective line $\hat{\mathbb{Q}}=\mathbb{Q}\cup\{\infty\}$ identify
$0$ with $0/1$ and $\infty$ with $1/0$. The set $\hat{\mathbb{Q}}$
forms the vertex set of $X$. Connect the (reduced) fractions $p/q$
and $r/s$ by an edge precisely when $|ps-rq|=1$. The resulting graph
is connected and each edge separates the graph.

The following result appears as Example 5.2 in \cite{Man}, where
Manning uses an interesting condition called the bottleneck
property.

\begin{prop}[\cite{Man}]\label{farey}
A Farey graph is quasi-isometric to a regular infinite-valence tree.
\end{prop}
Applying Example 1 from Section 2, we immediately obtain the
following result (see the remarks before Question 3 in Section 6 for
the lower bound).
\begin{cor} \label{cor2}The asymptotic dimension of a Farey
graph is $1$.\end{cor}

From Proposition \ref{farey} it is easy to see that a Farey graph is
not quasi-isometric to a finite valence tree. (See Proposition
\ref{fareynotproper} for a stronger statement.)  

\section{Property $A_1$}

G. Yu introduced property $A$ for discrete metric spaces as a means
to guarantee the existence of a uniform embedding into Hilbert space
\cite{Yu2}.  The existence of this embedding implies the coarse
Baum-Connes conjecture for discrete metric spaces with bounded
geometry. (As mentioned in the introduction many of the spaces we
consider are not examples of such spaces.) Yu's definition was the
following:

\begin{definition*} Let $\fin(Z)$ denote the collection of finite, nonempty
subsets of $Z$. The discrete metric space $Z$ has {\em property $A$}
if there are maps $A_n:Z\to \fin(Z\times\N)$ ($n=1,2,\ldots$) such
that
\begin{enumerate}
   \item for each $n$ there is some $R>0$ so that
\[A_n(z)\subset\{(z',j)\in Z\times\N\mid d(z,z')<R\}\]
for every $z\in Z$ and
   \item for every $K>0$
   \[\lim_{n\to\infty}\sup_{d(z,w)<K}\frac{|A_n(z)\Delta A_n(w)|}{|A_n(z)\cap A_n(w)|}=0.\]
\end{enumerate}
\end{definition*}

Following Tu \cite{Tu} and Dranishnikov \cite{Dr04}, define (for
$p\in\R\cup+\infty$) a metric space $X$ to have \emph{property
$A_p$} if there exist maps $a^n:X\to \ell^p(X)$ such that
$\|a^n_x\|_p=1,$ $a^n_x(y)\ge 0$ for all $x,y\in X$, $n\in\N,$ and
\begin{enumerate}
   \item there is a function $R=R(n)$ so that for all $x\in X,$
   $\supp(a^n_x)\subset N(x;R),$ and
   \item for every $K>0$,
   \[\lim_{n\to\infty}\sup_{d(x,y)<K}\|a^n_x-a^n_y\|_p=0.\]
\end{enumerate}

This definition is similar to a characterization of property $A$ for
discrete metric spaces with bounded geometry given by Higson and Roe
in \cite{HR}. There they show that for discrete metric spaces with
bounded geometry, property $A_1$ is equivalent to Yu's Property $A$.
Dranishnikov showed in \cite{Dr04} that discrete metric spaces with
property $A_p$ admit a coarse embedding in $\ell^p$.  In
\cite[Proposition 3.2]{Tu}, Tu shows that property $A_1$ and $A_2$
are equivalent. The proposition is stated for discrete metric
spaces, but the remarks following the proof state that the
equivalence holds in general. Also in remarks following this
proposition, Tu states that Property $A$ (in the sense of Yu)
implies property $A_1$. On the other hand, it is unknown whether
$A_1$ implies property $A$ without the requirement of bounded
geometry. In \cite[Proposition 3.2]{Dr04}, Dranishnikov shows that
properties $A_p$ and $A_q$ are equivalent for any $1\le p,q<\infty.$

The relation between finite asymptotic dimension and these
properties is the following:

\begin{lemma} {\rm\cite[Lemma 4.3]{HR}}\label{HR-Lemma} Let $X$ be a
discrete metric space with bounded geometry and finite asymptotic
dimension. Then $X$ has property $A$.
\end{lemma}

In fact, they show $X$ has property $A_1$, which is equivalent to
property $A$ for such spaces. Their argument along with Theorem
\ref{main} can be used to show the following result:

\begin{thm} Let $\G$ be a $\delta$-hyperbolic graph. Suppose $\sL$ is a set
of geodesics in $\G$ such that any two vertices of $\G$ are joined
by a geodesic in $\sL$. If $\sL$ satisfies property B with constants
$\ell,k=2 \delta$ and $D$, then $\G$ has property $A_1.$\end{thm}

We reproduce the argument from \cite{HR} for the reader's
convenience.  We also fill in some details using our definition of
asymptotic dimension in terms of $r$-multiplicity of the cover.

\begin{proof} By Theorem \ref{main}, $\as \G\le 2D-1.$  By the definition
of $\as,$ we know that for any $r\in\N$ we can find a uniformly
bounded cover $\sU$  of $\G$ with $5r$-multiplicity $\le 2D.$ Define
a cover $\sV$ of $\G$ by $\sV=\{N(U;2r)\mid U\in\sU\},$ where
$N(U;2r)$ denotes the $2r$-neighborhood of $U$ in $\G.$ Clearly
$\sV$ covers $\G$ and $\diam(V)\le 4r+\diam(U),$ so $\sV$ consists
of sets with uniformly bounded diameters. Next, if $x\in N(U;2r)$
then $N(x;5r)$ meets $U.$ For any $x$ there are at most $2D$ such
$U,$ so the multiplicity (or order) of $\sV$ is $\le 2D.$ Finally,
let $A\subset\G$ with $\diam(A)<r.$  If $x\in A,$ then $A\subset
N(U;2r)$ where $U\in\sU$ is any set containing $x.$ Thus, any set
with diameter $<r$ is entirely contained within some element of the
cover $\sV.$  (This says exactly that $r$ is a {\em Lebesgue number}
for $\sV.$)

For each $V\in\sV,$ define a map $\phi^r_V:\G\to [0,1]$ by
\[\phi^r_V(x)=\frac{d(x,\bar V)}{\sum_{W\in\sV}d(x,\bar W)},\]
where the bar denotes the complement in $\G$. Clearly, each
$\phi_V^r\not\equiv 0.$  For each $V\in\sV$ let $x_V$ denote some
point with $\phi_V^r(x_V)\neq 0.$  Now, define $a^r:\G\to
\ell^1(\G)$ by
\[a^r_x(z)=\sum_{V\in\sV}\phi_V^r(x)\delta_{x_{\empty_V}}(z),\]
where $\delta_{x_{\empty_V}}$ is the Dirac-$\delta$ function.
Clearly $a^r_x$ assumes a nonzero value at no more than $2D$ points
in $\G$ and so $a^r_x\in\ell^1(\G)$.

\noindent\textbf{Claim:}\begin{enumerate}
   \item[(1)] $a^r_x(z)\ge 0$ for all $z\in\G$;
   \item[(2)] $\|a^r_x\|_1=1$;
   \item[(3)] there is some $R=R(r)$ so that $\supp(a^r_x)\subset
   N(x;R)$ for all $x\in\G$; and
   \item[(4)] for every $K>0$ \[\lim_{r\to\infty}\sup_{d(z,w)<K}\|a^r_z-a^r_w\|_1=0.\]
\end{enumerate}

Item (1) is clear.  For item (2), we compute
\[\begin{array}{rcl}
\|a^r_x\|_1&=&\sum_{z\in\G}
\left|\sum_{V\in\sV}\phi^r_V(x)\delta_{x_{\empty_V}}(z)\right|\\
&=&\sum_{V\in\sV}\phi^r_V(x)=1.\\
\end{array}\]
To see item (3), we take $R(r)=4r+\diam(\sU)$ where $\diam(\sU)$ is
an upper bound for the diameters of the sets in $\sU$. Now,
$a^r_x(z)>0$ means that $z=x_V$ for some $V$ containing $x.$ Thus,
$z\in V$, so $d(z,x)\le R.$

Finally, for item (4), let $K>0$ and $\epsilon>0$ be given. Take $r$
so large that $(4D+1)^2K/r<\epsilon.$  The triangle inequality says
if $z,w\in\G$ and $V\in\sV,$ then $|d(z,\bar V)-d(w,\bar V)|\le
d(z,w).$ Also, observe that for any $w\in\Gamma,$
$\sum_{U\in\sV}d(w,\bar U)\ge r$ since $r$ is a Lebesgue number for
the cover $\sV.$ Thus, we have
\[|\phi_V^r(z)-\phi^r_V(w)|=
\left|\frac{d(z,\bar V)}{\sum_{U\in\sV}d(z,\bar U)}-\frac{d(w,\bar
V)}{\sum_{U\in\sV}d(w,\bar U)}\right|\]
\[\le \frac{\left|d(z,\bar V)-d(w,\bar V)\right|}{\sum_{U\in\sV}d(z,\bar U)}
+ \left|\frac{d(w,\bar V)}{\sum_{U\in\sV}d(z,\bar U)}-\frac{d(w,\bar
V)}{\sum_{U\in\sV}d(w,\bar U)}\right|\]
\[ \le \frac{d(z,w)}{\sum_{U\in\sV}d(z,\bar U)}+ \frac{d(w,\bar
V)}{\sum_{U\in\sV}d(z,\bar U)\sum_{U\in\sV}d(w,\bar
U)}\sum_{U\in\sV}\left|d(w,\bar U)-d(z,\bar U)\right|\]
\[\le\frac{1}{r}d(z,w)+
\frac{1}{r}\left(\sum_{U\in\sV}\left|d(z,\bar U)-d(w,\bar
U)\right|\right)\]\[\le \frac{1}{r}d(z,w)+\frac{4D
d(z,w)}{r}=\frac{(4D+1)}{r}d(z,w).\]

Then, \[ \begin{array}{rcl} \|a^r_z-a^r_w\|_1 & = &
\sum_{x\in\G}\left|a^r_z(x)-a^r_w(x)\right|\\
&=&\sum_{x\in\G}\left|
\sum_{V}\phi^r_V(z)\delta_{x_{\empty_V}}(x)-\sum_{V}\phi^r_V(w)
\delta_{x_{\empty_V}}(x)\right|\\
& = &
\left|\sum_V\left(\phi^r_V(z)-\phi^r_V(w)\right)\delta_{x_{\empty_V}}
\right|\\
& \le & \sup_{V}\left| \phi^r_V(z)-\phi^r_V(w)\right| 4D\\ & \le &
4D\frac{(4D+1)}r d(z,w)\le \frac{(4D+1)^2K}{r}< \epsilon.
\end{array}\]
So item (4) is proved. Thus, $\G$ has property $A_1.$
%
%
\end{proof}

As mentioned in Theorem 3, Masur-Minsky's tight geodesics $\sL$ on a
curve graph satisfy property B.  Thus we have the following result.

\begin{cor}
Let $S$ be an orientable surface of genus $g$ with $p$ punctures
such that $3g-3+p > 1$. Then the curve graph of $S$ has property
$A_1.$
\end{cor}

Tu \cite[Proposition 8.1]{Tu} also proved that a discrete hyperbolic
metric space $X$ with bounded geometry has property $A_1$ by fixing
$a\in\partial X$ and taking for the functions collections of
geodesic rays marching off to $a.$  As it stands, this argument
cannot be applied directly to discrete hyperbolic spaces with a
collection of geodesics satisfying property B as the following
example shows.

\begin{example*} Let $T$ be the following tree.  Let $x_0$ be a
vertex with infinite valence. Issuing from $x_0$ take edges of
length $1,2,3,\ldots.$  Clearly, $T$ is hyperbolic and the
collection of all geodesics in $T$ satisfies property B, but
$\partial T=\emptyset,$ so Tu's argument does not apply.
\end{example*}

On the other hand, suppose $\G$ is a discrete hyperbolic space with
a collection of geodesics $\sL$ with property B such that any two
points in $\G$ are connected by a geodesic in $\sL.$ If $\G$ has the
additional properties that $\partial\G\neq\emptyset$ and that
$\gamma\in\sL$ implies any subgeodesic of $\gamma$ is also in $\sL$,
then Tu's argument can be applied with very few changes to show that
$\G$ has property $A_1.$ The point of this remark is that when we
follow his argument we use only geodesics in the set $\sL$ to deal
with the problem that $\G$ may not have bounded geometry. We note
that since the set of tight geodesics of a curve graph satisfies
these additional properties, this approach gives another argument to
show that a curve graph has property $A_1$. After this paper was
written we were informed that Y.~Kida \cite{kida} has shown that a
curve graph has not only property $A_1$ but also property $A$.

A natural question at this point is whether the coarse Baum-Connes
conjecture is true for a curve complex. We settle this question by
showing that a curve graph is not coarsely equivalent to a proper
metric space, so the coarse Baum-Connes conjecture cannot be
formulated for it. It is likely that this is well known to
specialists, but we do not know a reference, so we record it.

The following lemma is well known.

\begin{lemma} Let $f:X\to Y$ be a coarse equivalence with coarse
inverse $g:Y\to X$. Then, for every $R$ there is an $S$ so that
$d(x,x')\ge S$ implies $d(f(x),f(x'))\ge R.$
\end{lemma}

\begin{proof} Since $X$ and $Y$ are coarsely equivalent, there is some
$K>0$ so that $d(fg,1_Y)\le K$ and $d(gf,1_X)\le K.$  Let $R$ be
given.  Since $g$ is a coarse map, there is a $S_g>0$ so that
$d(y,y')< R$ implies that $d(g(y),g(y'))< S_g.$  Put $S>2S_g+2K.$
Then, if $d(x,x')\ge S,$ we have $d(gf(x),gf(x'))\ge 2S_g+2K-2K.$
So, if $d(f(x),f(x'))< R,$ then we have $d(gf(x),gf(x'))< S_g$,
which is a contradiction. \end{proof}


Let $D>0$ be given. By an infinite $D$-discrete bounded subset of
$X$ we mean an infinite collection of points $x_i$ in $X$ such that
$d(x_i,x_j)\ge D$ whenever $i\neq j$ and $\sup\{d(x_i,x_j)\mid
i,j\}<\infty$.

\begin{prop}\label{prop S-non-proper}
Let $X$ be a metric space that has infinite $D$-discrete bounded
subsets for some $D$. Then $X$ is not proper.
\end{prop}

\begin{proof}
Take a sequence $\{x_i\}$ of distinct points in the infinite
$D$-discrete bounded subset. If $X$ were proper, this sequence would
have a convergent subsequence. Since $d(x_i,x_j)\ge D$ for all
$i\neq j$, this cannot happen.
\end{proof}

\begin{lemma} \label{lemma 3} Having an infinite $D$-discrete bounded subset for all
$D$ is an invariant of coarse isometry.
\end{lemma}
\begin{proof}
Let $f$ and $g$ be coarse maps from $X$ to $Y$ and $Y$ to $X$ that
give a coarse equivalence between $X$ and $Y$. Let $D>0$ be given.
Choose $D'>0$ so that if $d(x,x') \ge D'$  $(x,x' \in X)$ then
$d(f(x),f(x')) \ge D$. Let $A\subset X$ be an infinite $D'$-discrete
bounded subset. Then, $f(A)\subset Y$ is $D$-discrete by our choice
of $D'$ and bounded by the fact that $f$ is bornologous.
\end{proof}

\begin{cor}\label{fareynotproper}
Let $G$ be a Farey graph. Then $G$ is not coarsely equivalent to a
proper metric space.
\end{cor}

\begin{proof} Let $T$ be a tree with infinite valence.  We saw that $G$ and
$T$ are quasi-isometric. Since a coarse isometry between $G$ and a
proper metric space $X$ would yield a coarse isometry between $X$
and $T$, it suffices to show that $T$ and $X$ are not coarsely
isometric. So, we need only show that $T$ has infinite $D$-discrete
bounded subsets for any $D>0$.

To this end, let $D>0$ be given.
Fix a base vertex $t_0\in T$ and enumerate the infinitely many
vertices $t_1,t_2,\ldots$ adjacent to $t_0$. Let $\gamma_i$
($i=1,2,\ldots$) be geodesic rays from $t_0$ so that $\gamma_i$
contains $t_i$. Define $A=\{\gamma_i(D)\}$. Now, $A$ is an infinite
$2D$-discrete, $2D$-bounded subset of $T$. In particular, it is an
infinite $D$-discrete bounded subset of $T$.
\end{proof}

We remark that the above argument shows that there is no map $f:G
\to X$ such that there exist $C>0$ and $\epsilon \ge 0$ such that
for any $x,y \in G$, $\frac1C d(x,y) -\epsilon \le d(f(x),f(y)) \le
C d(x,y) + \epsilon$, (see the remarks following the definition of
coarse equivalence in Section 2).

\begin{lemma} \label{lemma 4}
$X=X(S_{g,p})$, $3g+p-4 > 0$, has infinite $D$-discrete bounded
subsets for all $D$.
\end{lemma}
\begin{proof}
We use the following theorem by Bowditch \cite{bowditch:tight}
concerning certain properness, which he calls {\it acylindricity},
of the action of $\Mod(S)$ on $X(S)$ (cf. the property WPD in
\cite{BF}; WPD is weaker, but enough for our purpose).

\begin{thm}{\rm\cite[Acylindricity]{bowditch:tight}}\label{acyl}
Let $X(S_{g,p})$, $3g+p-4 > 0$ be a curve graph. For any $D>0$,
there exist $L,K >0$, which depend on $D$ and $S_{g,p}$, with the
following property: if $x,y \in X$ are two points with $d(x,y) \ge
L$, then there are at most $K$ elements $a$ in $\Mod(S)$ such that
$d(x,ax) \le D$ and $d(y,ay) \le D$.
\end{thm}

Fix a surface $S$. Suppose $D>0$ is given. Let $L>0,K>0$ be
constants from Theorem \ref{acyl} for $D$ and $S$. Choose two points
$x,y \in X$ with $d(x,y) \ge L$. Let $G < \Mod(S)$ be the stabilizer
of $x$. $G$ is an infinite subgroup. For each $g \in G$, since
$d(gx,gy) \ge L$ there are at most $K$ elements $a \in G$ with
$d(gy,ay) \le D$ because of the acylindricity (note that
$d(gx,ax)=0$). It follows that, since $G$ is infinite, there are
infinitely many elements $g_i$ in $G$ such that $d(g_i y, g_j y) >D
(i \not= j)$. Since all points $g_i y $ are in $N(x;d(x,y)+1)$, we
have found the required subset.
\end{proof}

\begin{prop}\label{curvegraphnotproper}
A curve graph $X(S_{g,p})$, with $3g+p-4 > 0$ is not coarsely
equivalent to a proper metric space.
\end{prop}

\begin{proof} This is clear from Proposition \ref{prop S-non-proper},
Lemma \ref{lemma 3} and Lemma \ref{lemma 4}.
\end{proof}

Proposition \ref{curvegraphnotproper} was brought into our attention
by Misha Kapovich. His argument is different than the one we gave.


\section{Mapping class groups}
In this section we turn our attention to mapping class groups. The
reader is referred to Ivanov's paper \cite{Iv} for a thorough
introduction to mapping class groups. We reproduce many of the
definitions and results contained therein for the reader's
convenience. Let $S_{g,p}$ denote the compact orientable surface of
genus $g$ with $p$ punctures. The \emph{mapping class group of
$S_{g,p}$}, ${\rm Mod}(S_{g,p})$, is the group of isotopy classes of
orientation-preserving diffeomorphisms $S_{g,p}\to S_{g,p}.$ It is
also known by the name \emph{modular group}, which explains the
notation.

Let $\G$ be a finitely generated group. Fixing a (finite) set of
generators $S=S^{-1}$ endows $\G$ with a left-invariant word metric
defined by $d_S(g,h)=\|g^{-1}h\|_S$, where $\|g^{-1}h\|_S$ is the
length of the shortest $S$-word presenting the element $g^{-1}h.$
Notice that two finite generating sets give rise to quasi-isometric
metric spaces, so we define $\as \G$ to be the asymptotic dimension
of $(\G,d_S)$ where $S=S^{-1}$ is any finite generating set. 
We will have occasion to consider countable groups that are not
finitely generated. J.~Smith \cite{Smi} showed that every countable
group could be endowed with a left-invariant proper metric and that
all left-invariant proper metrics are coarsely equivalent. Thus,
when speaking of countable groups that may not be finitely
generated, we will always assume that the group has a proper
left-invariant metric on it. As a consequence, we see that $\as H\le
\as G$ for any subgroup $H$ of a finitely generated group $G$.

The mapping class group, ${\rm Mod}(S_{g,p})$ is finitely generated,
(in fact it is finitely presented, see \cite{Iv}). It was recently
shown that $\Mod(S)$ is exact \cite{Ha,Kida2}. As a consequence, the
Novikov conjecture holds for $\Mod(S)$ (cf. \cite{St}). A natural
question is whether these groups have finite asymptotic dimension
(see Question 2 in Section 7).

It is easy to obtain an obvious lower bound on $\as{\rm
Mod}(S_{g,p}),$ namely 
\[
\as{\rm Mod}(S_{g,p})\ge 3g-3+p.
\]
Indeed, there are $k=3g-3+p$ simple closed curves on $S_{g,p}$ and
Dehn twists by these curves commute because they are disjoint. Thus,
there is a copy of $\mathbb{Z}^{k}$ inside ${\rm Mod}(S_{g,p}).$
Dranishnikov, Keesling and Uspenskij showed that $\as
\mathbb{Z}^k=k$ in \cite{DKU}. Finally, if $Y\subset X$ in a metric
space $X,$ then any uniformly bounded cover of $X$ by sets with
$r$-multiplicity $\le n+1$ will restrict to a uniformly bounded
cover of $Y$ with $r$-multiplicity $\le n+1,$ so $\as Y\le \as X.$
Thus, we conclude that $\as{\rm Mod}(S_{g,p})\ge k.$

This lower bound is unlikely to give the exact asymptotic dimension
as we now demonstrate.

Observe that the Euler characteristic of $S_{g,p}$ is $2-2g-p.$ The
second part of \cite[Theorem 2.8.C]{Iv} states that when $2-2g-p<0,$
the following sequence is exact:
\[1\to \pi_1(S_{g,p})\to {\rm PMod}(S_{g,p+1})\to {\rm PMod}(S_{g,p})\to 1,\]
where ${\rm PMod}(S)$ is the \emph{pure mapping class group} of $S$,
defined as the group of all isotopy classes of all
orientation-preserving diffeomorphisms preserving setwise all
boundary components of $S$. 
Since ${\rm PMod}(S)$ is the kernel of the action of ${\rm Mod}(S)$
on the set of punctures, ${\rm PMod}(S)\subset{\rm Mod}(S)$ with
finite index.  Thus, $\as {\rm PMod}(S)=\as {\rm Mod}(S).$

A finitely presented group $\Gamma$ is said to be of type {\em FP}
if its classifying space $B\Gamma$ is homotopy dominated by a finite
complex. A finitely generated group $\Gamma$ is said to be of type
{\em VFP} if it contains an FP subgroup of finite index.

The following inequality holds for groups $\Gamma$ of type VFP,
\cite[Cor 4.11]{Dr05}:
\begin{equation}\label{dranish}{\rm vcd}\Gamma\le \as\G.\end{equation}
In \cite[\S5.4]{Iv}, Ivanov shows that $\Mod(S_{g,p})$ is of type
VFP when $2-2g-p<0.$


Ivanov computes the virtual cohomological dimension of mapping class
groups in \cite{Iv}, Theorems 6.4.A, 6.4.B and 6.4.C:
\begin{thm}[\cite{Iv}]\label{vcd} The following equalities hold for $S_{g,p}$,
the surface of genus g and p boundary components:
\[{\rm vcd}\Mod(S_{0,p})=\left\{
 \begin{array}{ll}
   0, & \hbox{if $p\le 3$;} \\
   p-3, & \hbox{if $p\ge 3$;} \\
\end{array}\right.\]
\[{\rm vcd}\Mod(S_{1,p})=\left\{
 \begin{array}{ll}
   1, & \hbox{if $p=0$;}\\
   p, & \hbox{if $p\ge 1$;}\\
\end{array}\right.\]
\[{\rm vcd}\Mod(S_{g,p})=\left\{
 \begin{array}{ll}
   4g-5, & \hbox{if $g\ge 2,\, p=0$;}\\
   4g-4+p, & \hbox{if $g\ge 2,\, p\ge 1$.}
 \end{array}
\right.\]
\end{thm}

Combining Theorem \ref{vcd} with Dranishnikov's result
\eqref{dranish} gives an improved lower bound on $\as\Mod(S_{g,p})$
in certain cases.

Next, we turn our attention to an upper bound for
$\as\Mod(S_{g,p}).$ Let $G$ be a finitely generated group and let
$1\to K\to G\to H\to 1$ be exact. Being the surjective image of $G$,
$H$ is finitely generated. On the other hand, $K$ is countable but
not necessarily finitely generated, so we consider any proper
left-invariant metric on $K$. For example, we could give $K$ the
metric it inherits as a subset of $G.$

Suppose $1\to K\to G\to H\to 1$ is exact and $G$ is finitely
generated. When $\as K$ and $\as H$ are both finite, Bell and
Dranishnikov prove that $\as G\le\as H+\as K$ in \cite{BD3}. If we
abuse notation slightly, allowing the terms of our inequalities to
be infinity, we get a formula
\[\as {\rm PMod}(S_{g,p+1})\le \as {\rm
PMod}(S_{g,p})+\as\pi_1(S_{g,p}),\] when $2-2g-p<0.$ Since $\as{\rm
PMod}(S)=\as{\rm Mod}(S),$ we get \[\as {\rm Mod}(S_{g,p+1})\le \as
{\rm Mod}(S_{g,p})+\as\pi_1(S_{g,p}),\] when $2-2g-p<0.$ Thus, we
can apply an inductive argument on the number of punctures of $S$.
We begin with an easy computation.

\begin{lemma} Let $S_{g,p}$ denote the compact surface with genus
$g$ and $p$ punctures with $g\ge 1$. Then
\[\as \pi_1(S_{g,p})=\left\{
                        \begin{array}{ll}
                          1, & \hbox{if $p>0$;} \\
                          2, & \hbox{if $p=0$.}
                        \end{array}
                      \right.\]
\end{lemma}

\begin{proof} Observe first that $\pi_1(S_{g,p})$ is a free group if $p>0,$
so in this case $\as\pi_1(S_{g,p})=1.$ If $p=0$ and $g>1,$ then
$\pi_1(S_{g,0})$ is quasi-isometric to $\mathbb{H}^2$, so
$\as\pi_1(S_{g,0})=2.$ Finally, $\pi_1(S_{1,0})=\mathbb{Z}\oplus
\mathbb{Z}$ so $\as\pi_1(S_{1,0})=2.$\end{proof}

\begin{thm} \label{Mod} Let $p\ge 0$ and $2-2g-p<0.$ If $\as\Mod(S_{g,0})<\infty,$ then
\[\as{\rm Mod}(S_{g,p})\le \as{\rm Mod}(S_{g,0})+p+1.\]
In particular, if $\as{\rm Mod}(S_{g,0})<\infty$ then $\as{\rm
Mod}(S_{g,p})<\infty$ for $p\ge 0.$
\end{thm}

\begin{proof} This follows easily from the lemma. Indeed,
\[
\begin{array}{ccl}
   \as {\rm Mod}(S_{g,p}) & \le & \as{\rm Mod}(S_{g,p-1}) + \as \pi_1(S_{g,p-1})\\
    & \vdots &  \\
    & \le & \as{\rm Mod}(S_{g,0})+ \as \pi_1(S_{g,p-1})+\cdots+\as \pi_1(S_{g,0})\\
     & = & \as{\rm Mod}(S_{g,0})+p+1
 \end{array}
\]\end{proof}

We can say more when $g\le 1.$  Ivanov explains in \cite[ 9.2]{Iv}
that ${\rm Mod}(S_{0,4})$ is commensurable with ${\rm
PSL}(2,\mathbb{Z}).$ Since ${\rm PSL}(2,\mathbb{Z})$ is
quasi-isometric to a tree, Example (i) in Section 3 implies that
$\as \Mod(S_{0,4})=1$. Ivanov also explains that ${\rm
Mod}(S_{1,1})\cong SL(2,\mathbb{Z})\cong
\mathbb{Z}_4\ast_{\mathbb{Z}_2}\mathbb{Z}_6$, and so $\as {\rm
Mod}(S_{1,1})=1$. Thus we obtain upper bounds for $\as{\rm
Mod}(S_{g,p})$ for $g=0$ or $1$. The lower bounds follow from
Theorem \ref{vcd} and Dranishnikov's estimate \eqref{dranish};
combining these bounds gives the following equalities.

\begin{cor}\label{g<2} If $p\ge 4$ then
\[
\as{\rm Mod}(S_{0,p})=p-3
\] and for $p\ge 1,$
\[
\as{\rm Mod}(S_{1,p})=p
\]
\end{cor}

Since the braid group $B_n$ is isomorphic to the mapping class group
of a disk with $n$ punctures, we see that a copy of $B_n$ sits
inside $\Mod(S_{0,n+1}),$ the mapping class group of the sphere with
$n+1$ punctures. Applying the previous result, we obtain the
following.

\begin{cor} \label{braid} Let $B_n$ be the braid group on $n$ strands.
If $n\ge 3,$ $\as B_n\le n-2$, so in particular $\as B_n<\infty$ and
$B_n$ has property A.
\end{cor}

Let $S$ be a finite set. A Coxeter matrix is a symmetric function
$M:S\times S\to\{1,2,3,\ldots\}\cup\{\infty\}$ with $m(s,s)=1$ and
$m(s,s')=m(s',s)\ge 2$ if $s\neq s'.$ The corresponding Coxeter
group $\sW(M)$ is the group with presentation
\[\sW(M)=\left<S\mid (ss')^{m(s,s')}=1\right>\] where $m(s,s')=\infty$
means no relation. The associated Artin group $\sA(M)$ is the group
with presentation
\[\sA(M)=\left<S\mid (ss')^{m(s,s')}=(s's)^{m(s,s')}\right>.\]
An Artin group with associated Coxeter matrix $M$, $\sA=\sA(M)$, is
said to be of finite type if $\sW(M)$ is finite. It is said to be of
affine type if $\sW(M)$ acts as a proper, cocompact group of
isometries on some Euclidean space with the elements of $S$ acting
as affine reflections.

In \cite{DJ}, it was shown that $\as\sW<\infty$ for a Coxeter group
$\sW.$ The following approach to finding an upper bound for the
$\as$ of (certain) Artin groups was suggested by Robert Bell.

In \cite{ChCr}, Charney and Crisp observe that each of the Artin
groups $\sA(A_n),$ $\sA(B_n)$ of finite type and the Artin groups
$\sA(\tilde{A}_{n-1})$ and $\sA(\tilde{C}_{n-1})$ of affine type is
a central extension of a finite index subgroup of $\Mod(S_{0,n+2})$
when $n\ge 3.$ Combining this with the fact that the centers of the
Artin groups of finite type are infinite cyclic and the centers of
those of affine type are trivial gives the following corollary.

\begin{cor}\label{artin-cor} Let $n\ge 3.$ Then if $\sA$ is
an Artin group of finite type $A_n$ or $B_n=C_n,$ we have $\as\sA\le
n$; if $\sA$ is an Artin group of affine type $\tilde{A}_{n-1}$ or
$\tilde{C}_{n-1},$ $\as\sA=n-1.$
\end{cor}

\begin{proof} By Corollary \ref{g<2}, $\as\Mod(S_{0,n+2})=n-1.$ Since their
centers are trivial, $\sA(\tilde{A}_{n-1})$ and
$\sA(\tilde{C}_{n-1})$ are themselves finite index subgroups of
$\Mod(S_{0,n+2}),$ we conclude that $\as\sA(\tilde{A}_{n-1})=n-1$
and $\as\sA(\tilde{C}_{n-1})=n-1$. For the Artin groups of finite
type the centers are infinite cyclic. Thus, $\sA(A_n)/Z$ and
$\sA(B_n)/Z$ are finite index subgroups of $\Mod(S_{0,n+2})$, so
$\as\sA(A_n)/Z=n-1,$ and $\as\sA(B_n)/Z=n-1.$ Since $\as$ of an
infinite cyclic group is 1, the extension theorem for $\as$ from
\cite{BD3} implies $\as \sA(A_n)\le n$ and $\as\sA(B_n)\le n.$
\end{proof}

\begin{cor} By Lemma \ref{HR-Lemma}, these groups have property
A.\end{cor}

%
%
%


Dan Margalit suggested the following to get an exact formula for the
$\as$ of the mapping class group of genus-2 surfaces. Let
$\Mod_2(S_{g,0})$ denote the elements of $\Mod(S_{g,0})$ that
commute with the hyperelliptic involution, (see
\cite{Big-Bud,Birman-Hilden}). Thinking of $\Mod_2(S_{g,0})\subset
\Mod(S_{g,0}),$ we see that $\as \Mod_2(S_{g,0})\le\as
\Mod(S_{g,0})$. In the genus-$2$ case, it is true that
$\Mod_2(S_{2,0})=\Mod(S_{2,0})$, (see \cite[Proposition
3.2]{Big-Bud} or \cite{Birman-Hilden}).

\begin{lemma} {\rm \cite{Big-Bud,Iv}} The following sequence is exact:
\[1\to\Z/2\Z\to\Mod_2(S_{2,0})\to \Mod(S_{0,6})\to 1.\]
\end{lemma}

\begin{cor}\label{cor 8} We have the following formula for the asymptotic dimension of
$\Mod(S_{2,p}):$
\[\as\Mod(S_{2,p})=\left\{
   \begin{array}{ll}
     3, & \hbox{if $p=0$;} \\
     p+4, & \hbox{if $p>0$.}
   \end{array}
 \right.\]
\end{cor}

\begin{proof} Obviously, $\as \Z/2\Z=0.$ By Corollary \ref{g<2}, we have
$\as\Mod(S_{0,6})=6-3=3.$ Applying the extension theorem for $\as$
from \cite{BD3} to the exact sequence in the previous lemma we see
that $\as\Mod_2(S_{2,0})\le 3.$ Thus, $\as\Mod(S_{2,0})\le 3.$
Finally, Theorem \ref{Mod} implies that $\as\Mod(S_{2,p})\le
p+1+3=p+4.$ The lower bounds follow from Theorem \ref{vcd} and
Dranishnikov's estimate \eqref{dranish}.

\end{proof}

%

It is interesting to observe that putting our lower bound 
for $\as {\rm Mod}(S_{g,p})$ together with Theorem \ref{Mod} we get
\[4g-4+p\le \as {\rm Mod}(S_{g,p})\le \as {\rm Mod}(S_{g,0})+1+p,\]
(when $2-2g-p<0$) so that if $\as{\rm Mod}(S_{g,0})$ is finite, then
$\as{\rm Mod}(S_{g,p})$ increases like the number of punctures.

\begin{lemma} \label{propA} Suppose $1\rightarrow K\xrightarrow{\iota}
G\xrightarrow{\phi} H\to 1$ is an exact sequence of groups with $G$
finitely generated. Suppose $\as H<\infty$ and $K$ has property A
with respect to some (hence any) proper, left-invariant metric. Then
$G$ has property $A$.
\end{lemma}

\begin{proof} Take a finite, symmetric generating set, $S=S^{-1}$, for $G$
and consider the word metric $d_S$ on $G.$ Taking the word metric
$d_{\phi(S)}$ corresponding to the generating set $\phi(S)$ on $H$
implies that $\phi$ is $1$-Lipschitz. Let $G$ act on $H$ by the rule
$g.h=\phi(g)h.$ In \cite{Be03}, the first author showed that if a
finitely generated group $G$ acts on a group $H$ with finite
asymptotic dimension so that for every $R,$ the set $W_R=\{g\in
G\mid d_{\phi(S)}(g.e,e)\le R\}$ has property $A$, then $G$ has
property $A$.

We will prove that $W_R=N(\iota(K);R),$ where the $R$-neighborhood
is taken in $G,$ with respect to the word metric in $G$. Thus, $W_R$
is quasi-isometric to $\iota(K).$ Thus $W_R$ has property $A$ if
$\iota(K)$ does. Since $K$ is countable, the metric space $\iota(K)$
and the group $K$ with its given proper metric are coarsely
equivalent. Since property $A$ is an invariant of coarse
equivalence, this will say that $W_R$ has property $A$ for every
$R,$ and hence,
that $G$ has property $A$. 

First assume $g\in W_R.$ Then $\|\phi(g)\|_{\phi(S)}\le R.$ So,
there exist $s_1,\ldots,s_k\in S,$ with $k\le R,$ so that
$\phi(g)=\phi(s_1)\cdots\phi(s_k).$ Take $g'=s_k^{-1}\cdots
s_1^{-1}\in G.$ Then, $gg'\in K,$ and $d(gg',g)\le R.$  On the other
hand, if $g\in N(K;R),$ then there is a $k\in K$ with $d(g,k)\le R.$
Since $\phi$ is $1$-Lipschitz, we have $d(\phi(g),e)\le R,$ as
required.

\end{proof}

The {\em Torelli group} $\I_g$ is the subgroup of ${\rm
Mod}(S_{g,0})$ acting trivially on $H_1(S_{g,0},\mathbb{Z})$, i.e.
the group $\I_g$ arising in the following exact sequence:\[1\to
\I_g\to {\rm Mod}(S_{g,0})\to Sp(2g,\mathbb{Z})\to 1.\]  Johnson
showed in \cite{Jo} that $\I_g$ is finitely generated if $g\ge 3,$
but $\I_2$ is non-finitely generated and free, \cite{Mess}. To make
sense of the asymptotic dimension of $\I_2,$ we endow it with a
proper left-invariant metric. One such metric on $\I_2$ is the one
it inherits as a subset of ${\rm Mod}(S_{2,0}).$ Since $\I_2$ is
isomorphic to a subgroup of $F_2$ the free group on two letters,
$\as\I_2\le 1.$ Since it contains a bi-infinite geodesic we can
conclude that $\as\I_2=1.$

\begin{cor} Endow the Torelli group $\I_g$ with a proper (left-invariant)
metric. Then $\I_g$ has finite asymptotic dimension if and only if
${\rm Mod}(S_{g,0})$ does. It has property $A$ if and only if ${\rm
Mod}(S_{g,0})$ does.
\end{cor}

\begin{proof} We showed in Corollary \ref{g<2} that $\as {\rm
Mod}(S_{g,0})<\infty$ when $g<2$ and so these groups (being finitely
generated) also have property $A.$ For $g<2$ the Torelli groups are
trivial
so they, too, have both properties.

If $g=2$ we only have to show that $\I_2$ has property $A$, but this
follows from the corresponding fact for the free group on two
generators, which contains it as a subgroup, and the fact that
property $A$ is a coarse invariant.

Let $g\ge 3.$ In both cases these properties pass to subsets. The
exact sequence
\[1\to \I_g\to {\rm Mod}(S_{g,0})\to Sp(2g,\mathbb{Z})\to 1\] and
the fact that $\as$ is finite for arithmetic groups by \cite{Ji}
implies the result for finite asymptotic dimension. Lemma
\ref{propA} implies the result for property $A$. \end{proof}

We briefly discuss Teichm\"uller spaces. Let $T(S)$ be the
Teichm\"uller space of the surface $S$ with negative Euler
characteristic.

See \cite{IT} or Chapter 5 in \cite{Iv} for definitions and
information. The group ${\rm Mod}(S)$ acts on $T(S)$ naturally. We
set $d(S_{g,p})=3g-3+p$. It is known that $T(S_{g,p})$ is
homeomorphic to ${\mathbb R}^m$, where $m=6g-6+2p=2d(S)$. There is a
natural metric on $T(S)$, called the {\it Teichm\"uller metric},
which is a Finsler metric. With respect to this metric, $T(S)$ is
proper, and the action of ${\rm Mod}(S)$ is proper and by
isometries. Therefore, we obtain
\begin{prop}\label{as-teich}
Suppose $S$ has negative Euler characteristic and endow $T(S)$ with
the Teichm\"uller metric. Then $\as {\rm Mod}(S) \le \as T(S)$,
where we allow the possibility that the terms of this inequality are
infinite.
\end{prop}
\begin{proof} The group $\Mod(S)$ acts on $T(S)$ (with the Teichm\"uller
metric) properly by isometries \cite[Chapter 6.3]{IT}, so
Proposition 2.3 from \cite{Ji} applies to this setting. We outline
an argument. We fix a point $x \in T(S)$ and identify each element
$g \in \Mod(S)$ with $g.x \in T(S)$. We wish to consider the metric
on $\Mod(S)$ induced from being a subset in $T(S)$, but this is a
pseudo-metric as we could have non-trivial $g$ with $g.x=x.$ Since
the action is proper, such a $g$ must be torsion, so we may (if
necessary) pass to a torsion-free subgroup of $\Mod(S)$ of finite
index, \cite{Iv}. This subgroup has the same $\as$ as $\Mod(S).$ The
point is that the metric this group inherits as a subset of $T(S)$
is a coarsely equivalent to a word metric on $\Mod(S).$ Thus, the
asymptotic dimension of $\Mod(S)$ as a finitely generated group is
the same as the asymptotic dimension of $\Mod(S)$ as a subset of
$T(S).$ It follows that $\as {\rm Mod}(S) \le \as T(S)$.\end{proof}

\section{Some open problems}

Our results lead to some natural questions. Again, we let $S_{g,p}$
be the compact orientable surface of genus $g$ with $p$ punctures
and let $X(S_{g,p})$ denote its curve graph.
\begin{question}
What is the asymptotic dimension of $X(S_{g,p})$?
\end{question}
It is known that $\dim(C(S_{g,p}))=3g+p-4$ (in the non-exceptional
cases where this number is positive). Unfortunately, generally there
is no relation between $\dim$ and $\as$; for instance, $\as C(S)=\as
C^{(1)}(S)$. We showed that $\as {\rm Mod}(S_{g,p})$ (if finite)
depends on $g$ and $p$, (see the remarks following Corollary
\ref{cor 8}). It would be very interesting if $\as X(S_{g,p})$ were
independent of $g$ and $p$.

Perhaps a promising approach to computing the asymptotic dimension
of a curve graph would be to examine the boundary $\partial
X(S_{g,p})$. For finitely generated hyperbolic groups $\G$, it is
known that $\dim\partial\G+1=\as\G$ by results of Buyalo and
Lebedeva \cite{BL}, cf. \cite{BuSch,Swi,Gr93}. In fact, Buyalo and
Lebedeva \cite{BL} prove this formula not only for groups, but for
cobounded, hyperbolic, proper, geodesic metric spaces. The drawback
to this approach is that the boundary of $X(S_{g,p})$ is not well
understood.

If $3g+p-3>1$, $X(S_{g,p})$ contains a bi-infinite geodesic, (cf.
\cite{bowditch:tight}) so $\as X(S_{g,p})\ge 1$. Note that the
asymptotic dimension of a Farey graph is $1$, (Corollary
\ref{cor2}).
%

Although we had some finiteness results for $\as{\rm Mod}(S_{g,p})$
the question remains open when $g\ge 3.$
\begin{question}\label{mod-finite}
Is the asymptotic dimension of ${\rm Mod}(S_{g,p})$ finite when
$g\ge 3$?
\end{question}

Notice that if $\as\Mod(S)<\infty$ then by Higson and Roe's result
\cite[Lemma 4.3]{HR}, $\Mod(S)$ is exact. Thus, the result on the
exactness of $\Mod(S)$ \cite{Ha,Kida2} would follow from an
affirmative answer to this question.

A naive approach to this question would be to attempt to use the
Hurewicz-type theorem of Bell and Dranishnikov \cite{Be03,BD3}. The
group ${\rm Mod}(S_{g,p})$ acts on $X(S_{g,p})$ by isometries. Since
$X(S_{g,p})$ has finite asymptotic dimension, we would be able to
conclude that ${\rm Mod}(S_{g,p})$ has finite asymptotic dimension
(or property $A$) provided we could show that the set $\{g\in {\rm
Mod}(S_{g,p})\mid d_X(g.x_0,x_0)\le r\}$ has finite asymptotic
dimension (respectively, property $A$) for all $r\in\N.$ Here, $x_0$
is any point of $X$. We remark that if we regard $x_0$ as a curve on
$S$ the set we need to analyze contains the stabilizer subgroup of
$x_0$, which is ${\rm Mod}(S \backslash x_0)$. A problem is that
those two sets are not necessarily coarsely equivalent (cf. proof of
Lemma \ref{propA}).

L.~Ji proved, in \cite{Ji}, that arithmetic groups have finite
asymptotic dimension. Ivanov draws a comparison between mapping
class groups and arithmetic groups in Chapter 9 of \cite{Iv}.
Although the groups are different (cf. \cite{BF}), many things that
are proved true for arithmetic groups are later proved to be true
for mapping class groups. So, an affirmative answer to Question
\ref{mod-finite} does not seem unlikely. Ji gave an upper bound of
the asymptotic dimension of an arithmetic group, $\Gamma$, as
follows (cf. proof of Proposition \ref{as-teich}). Let $X$ be the
symmetric space associated to $\Gamma$. Then $\as \Gamma \le \as X
=\dim X$.
So, in connection with Question \ref{mod-finite}, we ask the
following (see Proposition \ref{as-teich}).
\begin{question}\label{teich-finite}
Let $S=S_{g,p}$ with negative Euler characteristic and endow $T(S)$
with the Teichm\"uller metric. What is $\as T(S)$? Is it finite? Is
$\as T(S)$ strictly bigger than $\as \Mod(S)$? If not (i.e., if they
are same), then is $\Mod(S)$ coarsely equivalent to $T(S)$ ?
\end{question}

We can ask similar questions on the asymptotic dimensions of a
Teichm\"uller space with the Weil-Petersson metric (cf.\cite{IT}),
$g_{\rm WP}$, but we do not even know if a statement similar to
Proposition \ref{as-teich} holds because the space is not proper.
Let $n$ be the greatest integer less than or equal to $(d(S)+1)/2.$
Brock and Farb (Theorem 1.4 \cite{BrFa}) proved that there is a
quasi-isometric embedding $\mathbb{R}^n \to (T(S_{g,p}), g_{\rm
WP}).$ It follows that $\as
(T(S), g_{\rm WP}) \ge n.$ 

Let $S, S'$ be compact orientable surfaces. If there exists a system
of disjoint curves, $C$, on $S$ so that one of the connected
components of $S \backslash C$ is homeomorphic to the interior of
$S'$, we write $S' < S$. Then we have $\as (T(S'), g_{\rm WP}) \le
\as (T(S), g_{\rm WP})$. This follows from the fact that $(T(S),
g_{\rm WP})$ contains a subset which is quasi-isometric to $(T(S'),
g_{\rm WP})$. Indeed (cf. \cite{Wo}), although the metric space
$(T(S), g_{\rm WP})$ is not complete, its metric completion
$\overline{(T(S), g_{\rm WP})}$ contains $(T(S'), g_{\rm WP})$ as a
totally geodesic subspace in such a way that $(T(S'), g_{\rm WP})
\subset \overline{(T(S), g_{\rm WP})}\backslash (T(S), g_{\rm WP})$.
For $\epsilon
>0$, let $N$ be the $\epsilon$-neighborhood of $T(S')$ in
$\overline{T(S)}$. Then $(N \backslash T(S')) \subset T(S)$ is
quasi-isometric to $T(S')$.

\end{document}